# A Triangular Array of the Counts of Natural Numbers with the Same Number of Prime Factors (Dimensions) Within $2^n$ Space

**Abstract**

By defining the dimension of natural numbers as the number of prime factors, all natural numbers smaller than $2^{n+1}$ ($n \in \mathbf{N}$) can be classified by their dimensions, and the count of numbers of each dimension gives a 'dimensions distribution' directly related to the distribution of prime numbers. A triangle similar to Pascal's triangle from some aspects can be obtained by the ensemble of 'dimensions distributions'. Some superficial explorations have been done to this interesting triangle.

## 1. The dimension of natural numbers

Every natural number > 1 can be expressed as the product of its prime factor(s) and 1. If we construct a coordination system with only prime numbers on the axis as scales and set 1 to be the origin (Fig 1), then natural numbers with $N$ prime factors can be presented in the corresponding $N$-dimensional coordination system. In this way we may define the 'dimension' of a natural number as the number of its prime factors. Subsequently prime numbers are 1D number and 1 is defined as zero-dimensional number. The natural number 3, 10, 18 as examples of 1D, 2D, 3D number and their presentations in the 3D prime coordination system are shown in Fig 1.

## 2. Dimensions distribution within $2^n$ space

Obviously along the number line the first $n$-dimensional number is $2^n$, then all natural numbers $\leq 2^n$ must be able to be presented in the $n$-dimensional coordination. This coordination system is confined by the size of $2^n$, therefore we name the natural number range $[1, 2^n]$ as $2^n$ space. The count of numbers of each dimensions $\leq n$ assemble a distribution fulfilling the $2^n$ space. In any $2^n$ space with $n \geq 1$, we have and only have one 0D number (1) and one ($n$-th)D number ($2^n$). For example, for the space $2^5$ we have:
one 0D number 1;
eleven 1D numbers (Prime) 2, 3, 5, 7, 11, 13, 17, 19, 23, 29, 31;
ten 2D numbers 4, 6, 9, 10, 14, 15, 21, 22, 25, 26;
seven 3D numbers 8, 12, 18, 20, 27, 28, 30;
two 4D numbers 16, 24;
one 5D number $2^5=32$. The distribution of dimensions within the $2^5$ space is shown in Fig 2.

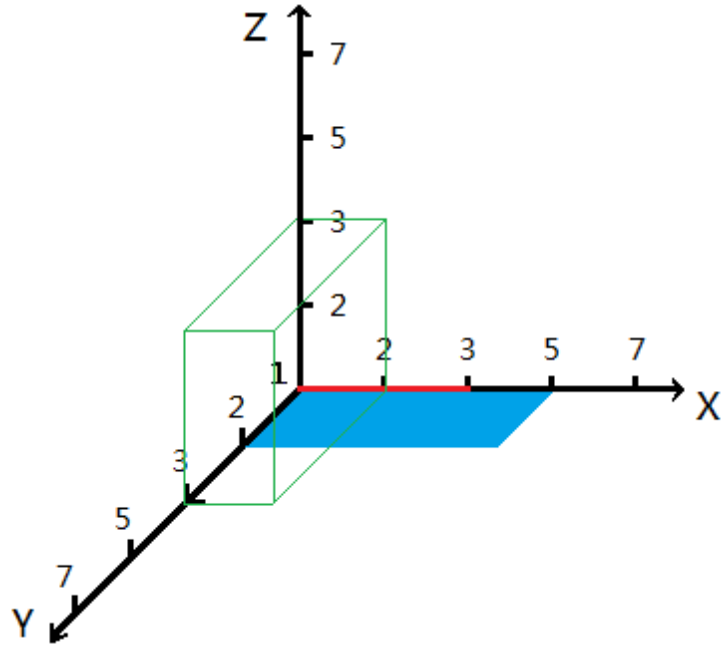

Fig 1. A 3D coordination system to present natural numbers with only prime numbers on the axis as scales and 1 to be the origin. The red line indicates the number 3, which is one dimensional on the X axis. The navy-blue area indicates 10, which is 2D presented on the X-Y plate. The green lined box indicates 18 and it is a 3D number presented in the X-Y-Z space.

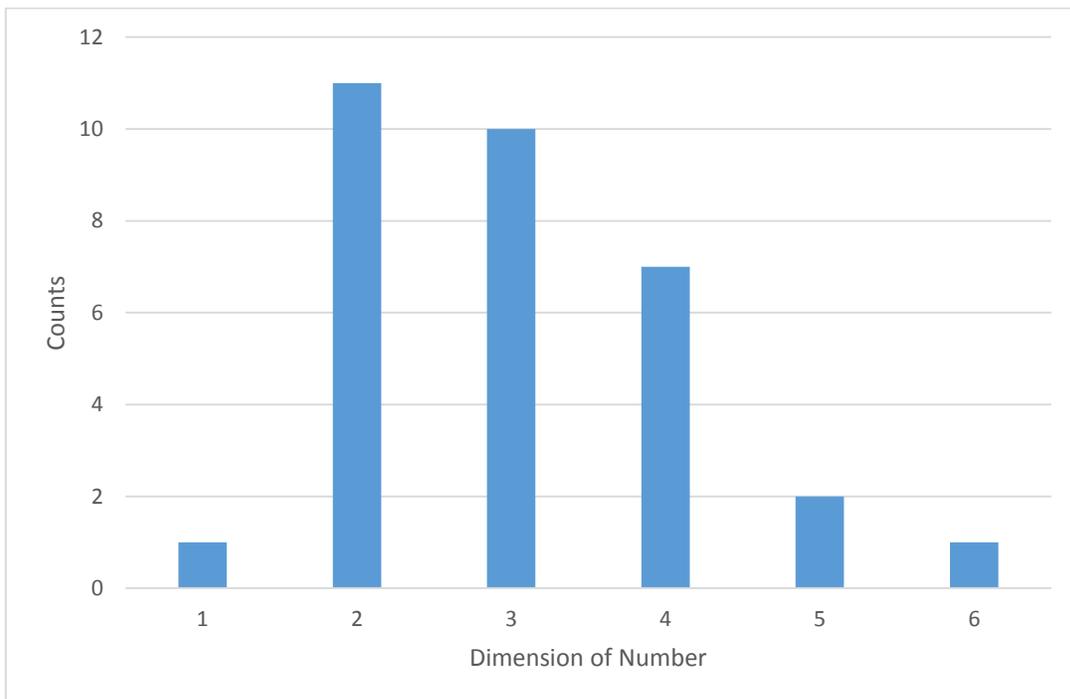

Fig 2. The distribution of dimensions within the $2^5$ space

## 3. The triangle of dimensions distribution

Although for any natural number range [1, $n$] ($n \in \mathbf{N}$) we can have a distribution of dimensions, here we take $2^n$ spaces as our interests. By accounting the dimensions distributions of $2^n$ spaces as $n$ = 1, 2, 3 … and so on, we can fit all the distributions in a triangle which has dimensions as columns and $2^n$ spaces as rows. The triangle within the $2^{23}$ space is shown in table 1.

Several interesting phenomena in this triangle:
1) It counts all the natural numbers, the sum of terms in $m$-th row equals to $2^m$, this property is as the same as the Pascal's triangle [1].
2) For each row, the last several terms are fixed, for example for any $2^n$ ($n > 5$) space, the $n$-th column is 1, the ($n$-1)-th column is 2 and the ($n$-2)-th column is 7.
3) The triangle of size $2^{16}$ of all the distributions is shown as a 3D map in Fig 3, since the number is much larger within higher spaces, for a better observation the logarithm of all the counts is graphed in Fig 4 (3D map in Fig 4A and 2D graph in Fig 4B).

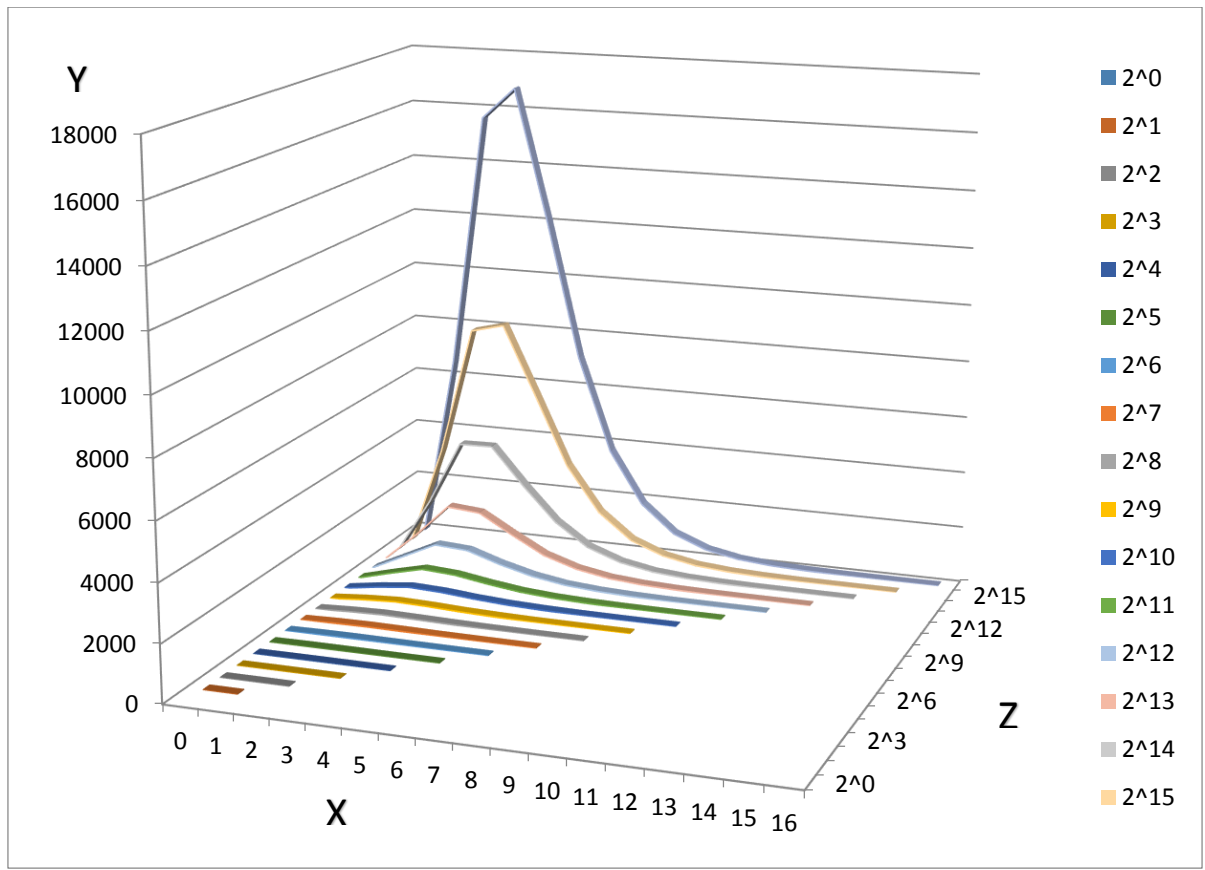

Fig 3. The triangle of size $2^{16}$ of all the distributions shown as 3D map. The X axis is dimensions, Y axis is the counts of numbers and the distribution of each $2^n$ spaces are arranged on Z axis.

Table 1. The triangle of Dimensions Distribution within $2^{23}$

| | 0 | 1 | 2 | 3 | 4 | 5 | 6 | 7 | 8 | 9 | 10 | 11 | 12 | 13 | 14 | 15 | 16 | 17 | 18 | 19 | 20 | 21 | 22 | 23 |
|---|---|---|---|---|---|---|---|---|---|---|---|---|---|---|---|---|---|---|---|---|---|---|---|---|
| 2^0 | 1 | | | | | | | | | | | | | | | | | | | | | | | |
| 2^1 | 1 | 1 | | | | | | | | | | | | | | | | | | | | | | |
| 2^2 | 1 | 2 | 1 | | | | | | | | | | | | | | | | | | | | | |
| 2^3 | 1 | 4 | 2 | 1 | | | | | | | | | | | | | | | | | | | | |
| 2^4 | 1 | 6 | 6 | 2 | 1 | | | | | | | | | | | | | | | | | | | |
| 2^5 | 1 | 11 | 10 | 7 | 2 | 1 | | | | | | | | | | | | | | | | | | |
| 2^6 | 1 | 18 | 22 | 13 | 7 | 2 | 1 | | | | | | | | | | | | | | | | | |
| 2^7 | 1 | 31 | 42 | 30 | 14 | 7 | 2 | 1 | | | | | | | | | | | | | | | | |
| 2^8 | 1 | 54 | 82 | 60 | 34 | 15 | 7 | 2 | 1 | | | | | | | | | | | | | | | |
| 2^9 | 1 | 97 | 157 | 125 | 71 | 36 | 15 | 7 | 2 | 1 | | | | | | | | | | | | | | |
| 2^10 | 1 | 172 | 304 | 256 | 152 | 77 | 37 | 15 | 7 | 2 | 1 | | | | | | | | | | | | | |
| 2^11 | 1 | 309 | 589 | 513 | 325 | 168 | 81 | 37 | 15 | 7 | 2 | 1 | | | | | | | | | | | | |
| 2^12 | 1 | 564 | 1124 | 1049 | 669 | 367 | 177 | 83 | 37 | 15 | 7 | 2 | 1 | | | | | | | | | | | |
| 2^13 | 1 | 1028 | 2186 | 2082 | 1405 | 770 | 392 | 182 | 84 | 37 | 15 | 7 | 2 | 1 | | | | | | | | | | |
| 2^14 | 1 | 1900 | 4192 | 4214 | 2866 | 1643 | 831 | 406 | 185 | 84 | 37 | 15 | 7 | 2 | 1 | | | | | | | | | |
| 2^15 | 1 | 3512 | 8110 | 8401 | 5931 | 3410 | 1790 | 867 | 414 | 186 | 84 | 37 | 15 | 7 | 2 | 1 | | | | | | | | |
| 2^16 | 1 | 6542 | 15658 | 16771 | 12139 | 7150 | 3757 | 1880 | 887 | 418 | 187 | 84 | 37 | 15 | 7 | 2 | 1 | | | | | | | |
| 2^17 | 1 | 12251 | 30253 | 33427 | 24782 | 14859 | 7942 | 3973 | 1931 | 900 | 420 | 187 | 84 | 37 | 15 | 7 | 2 | 1 | | | | | | |
| 2^18 | 1 | 23000 | 58546 | 66550 | 50444 | 30769 | 16677 | 8432 | 4102 | 1962 | 907 | 421 | 187 | 84 | 37 | 15 | 7 | 2 | 1 | | | | | |
| 2^19 | 1 | 43390 | 113307 | 132405 | 102458 | 63528 | 34820 | 17820 | 8736 | 4179 | 1979 | 911 | 421 | 187 | 84 | 37 | 15 | 7 | 2 | 1 | | | | |
| 2^20 | 1 | 82025 | 219759 | 262865 | 207945 | 130713 | 72515 | 37423 | 18532 | 8917 | 4225 | 1989 | 913 | 421 | 187 | 84 | 37 | 15 | 7 | 2 | 1 | | | |
| 2^21 | 1 | 155611 | 426180 | 522296 | 420511 | 268776 | 150454 | 78369 | 39045 | 18968 | 9028 | 4250 | 1995 | 914 | 421 | 187 | 84 | 37 | 15 | 7 | 2 | 1 | | |
| 2^22 | 1 | 295947 | 827702 | 1036033 | 850518 | 550342 | 311933 | 163456 | 82065 | 40051 | 19230 | 9094 | 4266 | 1998 | 914 | 421 | 187 | 84 | 37 | 15 | 7 | 2 | 1 | |
| 2^23 | 1 | 564163 | 1608668 | 2055256 | 1716168 | 1126507 | 643941 | 340656 | 171749 | 84366 | 40670 | 19388 | 9133 | 4274 | 2000 | 914 | 421 | 187 | 84 | 37 | 15 | 7 | 2 | 1 |

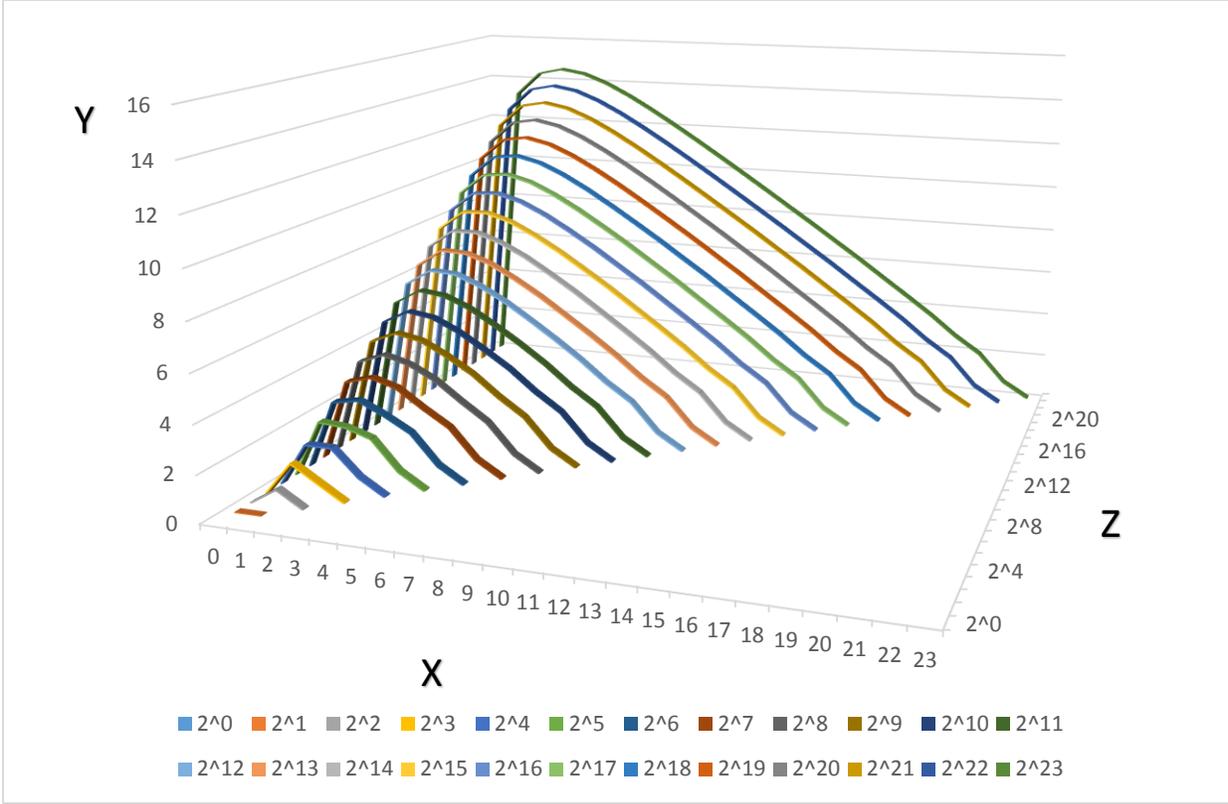

(A)

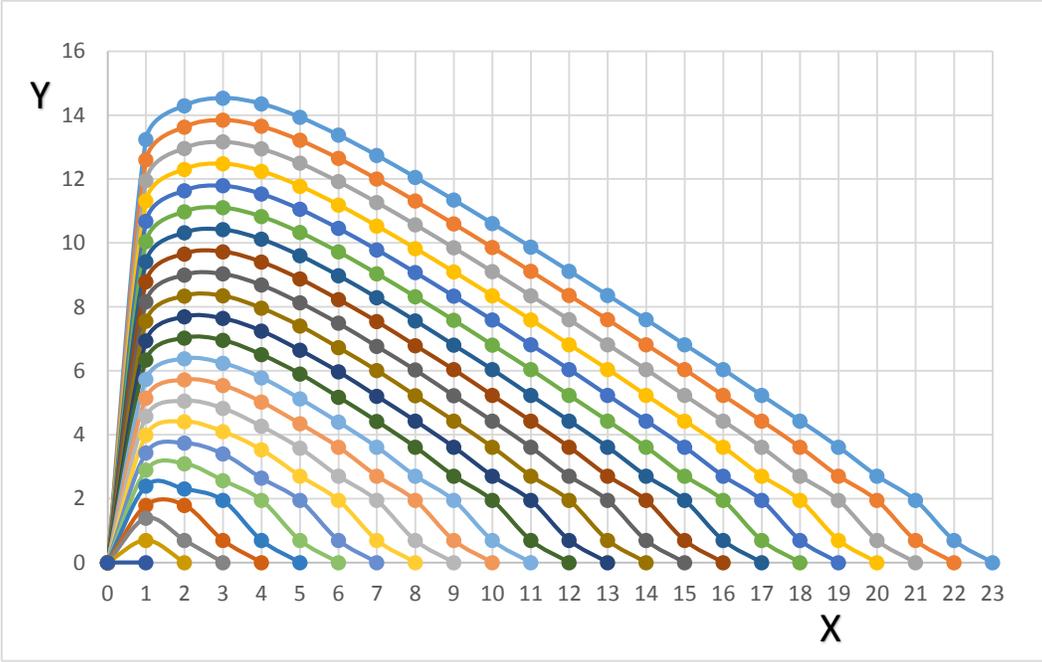

(B)

Fig 4. The logarithm of all the counts of $2^{16}$ triangle shown in 3D (A) and 2D (B)

4) For $n \leq 5$, the prime column in $2^n$ space has the most numbers, that is, the maximum column. For $6 \leq n \leq 13$, the 2D is the maximum column in $2^n$ space. For $14 \leq n \leq 23$, the 3D is the maximum column. Further calculations with $n > 23$ were not done.

Question 1: Will the maximum column continue to shift to right or fixed on the 3D column?

## 4. Analysis of the phenomenon 3-2)

It is easy to understand the fixed last several columns of each row: the fixed $(n-x)$-th $(x \in \mathbf{N})$ term in the $2^n$ space (the $(n+1)$-th row in the triangle) should be determined by the number of primes in the $2^{x+1}$ space. Take $x = 2$ as an example, in Fig 5 we draw a diagonal arrow starting from (row: $2^3$, column: 1) and crossing each $(n-2)$-th column of all following $2^n$ spaces.

|       | 0 | 1  | 2  | 3  | 4  | 5  | 6 | 7 | 8 |
|-------|---|----|----|----|----|----|---|---|---|
| 2^0   | 1 |    |    |    |    |    |   |   |   |
| 2^1   | 1 | 1  |    |    |    |    |   |   |   |
| 2^2   | 1 | 2  | 1  |    |    |    |   |   |   |
| 2^3   | 1 | 4  | 2  | 1  |    |    |   |   |   |
| 2^4   | 1 | 6  | 6  | 2  | 1  |    |   |   |   |
| 2^5   | 1 | 11 | 10 | 7  | 2  | 1  |   |   |   |
| 2^6   | 1 | 18 | 22 | 13 | 7  | 2  | 1 |   |   |
| 2^7   | 1 | 31 | 42 | 30 | 14 | 7  | 2 | 1 |   |
| 2^8   | 1 | 54 | 82 | 60 | 34 | 15 | 7 | 2 | 1 |
| ...   |   |    |    |    |    |    |   |   |   |

Fig 5. The diagonal arrow starting from (row: $2^3$, column: 1) and crossing each $(n-2)$-th column of all following $2^n$ spaces

The numbers appear on the diagonal line are: 4, 6, 7, 7… that is, in any $2^n$ ($n \geq 5$) space, the number of $(n-2)$-dimensional numbers are fixed to be seven. In the space of $2^3$, the four prime numbers are: 2, 3, 5, 7. Now we count the number of $(n-2)$-th dimensional numbers in a $2^n$ ($n \geq 5$) space:

1) In this space the first $(n-2)$D number must be $2^{n-2}$.

2) The second $(n-2)$D number is $2^{n-3} \times 3$. The third one is $2^{n-3} \times 5$, the fourth one is $2^{n-3} \times 7$. And we don't have as the fifth one as $2^{n-3} \times 11$, because $2^{n-3} \times 11$ is larger than the $2^n = 2^{n-3} \times 2^3$ space, or simpler to say, 11 is larger than $2^3 = 8$. It is clear that for the combination form of ($2^{n-3}$ times Prime), the Prime must be smaller than $2^3$ (which is from $2^n / 2^{n-3}$). Note that although 2 is also a Prime number, $2^{n-3} \times 2$ has been accounted in the previous step.

3) Now taking $2^{n-4}$ to make the combinations, we have $2^{n-4} \times 3 \times 3$, $2^{n-4} \times 3 \times 5$, as the fifth and sixth $(n-2)$D number, no more combination form of ($2^{n-4}$ times 2D) is available, because similarly the 2D number here is limited by the size of $2^4$ (which is from $2^n / 2^{n-4}$). Note that $2^{n-4} \times 2 \times 2$ was accounted in 4-1), and $2^{n-4} \times 2 \times 3$, $2^{n-4} \times 2 \times 5$, $2^{n-4} \times 2 \times 7$ were accounted in 4-2).

4) For the combination form of ($2^{n-5}$ times 3D), we only have $2^{n-5} \times 3 \times 3 \times 3$ as the seventh $(n-2)$D number, $3 \times 3 \times 5$ and all other non-2-prime-factor 3D numbers are out of the limit of $2^5$.

5) No $(n-2)$D number is available in the form of ($2^{n-6}$ times 4D), $3^4$ is larger than $2^6$.

Therefore in any $2^n$ ($n \geq 5$) space, the number of $(n-2)$-dimensional numbers are fixed to be seven. The step 4-3) and 4) also reveals the reason of the number six in the space $2^4$ on the arrow, because the $2^4$ space is impossible to have the combination form of ($2^{n-5}$ times 3D), the seventh number accounted in 4-4) is ruled out.

Basically the triangle is a picture of prime numbers and prime factors, and it directly relates to the prime number distribution. It is clear that for the similar arrows on this triangle:
4, 6, 7, 7…
6, 10, 13, 14, 15, 15…
11, 22, 30, 34, 36, 37, 37…
and so on, the first number (number of primes in $2^{x+1}$ space) determines the rest terms in the series. The terms following the number of primes also depend on the size of the space, and at some point when the space is large enough, the eventual term will be fixed regardless of the space size.

We may hypothesis that the terms in each series would be a function of the first term and the space size, however the author is not able to figure out this function. If this function actually exists in a solvable form, the exact distribution of prime number in $2^n$ is easy to approach by calculating the number of composite numbers (non-1D numbers) from the information of $2^{n-1}$ space (previous row), the calculation process is shown in Fig 6.

Question 2: Do this function of diagonal series and its subsequent algebra of prime distribution exist?

## 5. Analysis of phenomenon 3-3)

The Fig 4B shows a quite well-shaped map, for each $x \geq 2$ (column in the triangle), the logarithm of counts (Y value) seems increased principally. Table 2 shows the data set of Y value on $x = 1, 2, 3, 4, 10$ and $13$ in Fig 4B. The difference of neighbor Y values are calculated next to the corresponding column and mapped in Fig 7A, B, C, D, E, and F.

|      | 0 | 1  | 2  | 3  | 4  | 5  | 6 | 7 | 8 |
|------|---|----|----|----|----|----|---|---|---|
| 2^0  | 1 |    |    |    |    |    |   |   |   |
| 2^1  | 1 | 1  |    |    |    |    |   |   |   |
| 2^2  | 1 | 2  | 1  |    |    |    |   |   |   |
| 2^3  | 1 | 4  | 2  | 1  |    |    |   |   |   |
| 2^4  | 1 | 6  | 6  | 2  | 1  |    |   |   |   |
| 2^5  | 1 | 11 | 10 | 7  | 2  | 1  |   |   |   |
| 2^6  | 1 | 18 | 22 | 13 | 7  | 2  | 1 |   |   |
| 2^7  | 1 | 31 | 42 | 30 | 14 | 7  | 2 | 1 |   |
| 2^8  | 1 | 54 | 82 | 60 | 34 | 15 | 7 | 2 | 1 |

Fig 6. Calculation of prime distribution in $2^n$ spaces by the hypothesized function of diagonal series. The red arrow represents the hypothesis function calculation. The blue arrow shows that once all the number of composite numbers can be calculated from the previous rows, the number of primes can be obtained by subtraction.

It can be observed in Fig 7 that, except the abnormal behavior of primes, other dimensional numbers approximately obey the rule that after an initial oscillation the increase of Y value approaches to a constant, especially for small dimensions.

Question 3: For each composite number column in Fig 4B, does there exist a (or approximate) constant increase to make the column to be a disciplinary series? Or is there a more general principle for the variation of distributions?

## 6. Entropy of dimensions distribution

A main feature of primes is that the appearance of prime numbers on number line looks like random, however the distribution of prime numbers seems to obey some principle and the research on Prime Number Theorem is still in progress.[2] Since this triangle is basically constructed on the primes, we would like to quantify the randomness of the dimensions distribution. The entropy of each dimension distribution within $2^{23}$ space is calculated by the Gibbs definition in thermodynamics[3]:
$S = -k_B \sum P_i \text{Ln} P_i$, the constant $k_B$ is set to be 1.

The entropies of the distribution of each row in Pascal's triangle are also shown in Fig 8 for comparison.

Table 2. The data set of Y value on *x* = 1, 2, 3, 4, 10 and 13 in Fig 4B

|      | 1        | difference | 2        | difference | 3        | difference |
|------|----------|------------|----------|------------|----------|------------|
| 2^0  |          |            |          |            |          |            |
| 2^1  | 0        |            |          |            |          |            |
| 2^2  | 0.693147 |            | 0        |            |          |            |
| 2^3  | 1.386294 | 0.693147   | 0.693147 |            | 0        |            |
| 2^4  | 1.791759 | 0.405465   | 1.791759 | 1.098612   | 0.693147 |            |
| 2^5  | 2.397895 | 0.606136   | 2.302585 | 0.510826   | 1.94591  | 1.252763   |
| 2^6  | 2.890372 | 0.492476   | 3.091042 | 0.788457   | 2.564949 | 0.619039   |
| 2^7  | 3.433987 | 0.543615   | 3.73767  | 0.646627   | 3.401197 | 0.836248   |
| 2^8  | 3.988984 | 0.554997   | 4.406719 | 0.669050   | 4.094345 | 0.693147   |
| 2^9  | 4.574711 | 0.585727   | 5.056246 | 0.649527   | 4.828314 | 0.733969   |
| 2^10 | 5.147494 | 0.572783   | 5.717028 | 0.660782   | 5.545177 | 0.716864   |
| 2^11 | 5.733341 | 0.585847   | 6.378426 | 0.661398   | 6.240276 | 0.695098   |
| 2^12 | 6.335054 | 0.601713   | 7.024649 | 0.646223   | 6.955593 | 0.715317   |
| 2^13 | 6.935370 | 0.600316   | 7.689829 | 0.665180   | 7.641084 | 0.685492   |
| 2^14 | 7.549609 | 0.614239   | 8.340933 | 0.651105   | 8.346168 | 0.705083   |
| 2^15 | 8.163941 | 0.614332   | 9.000853 | 0.659920   | 9.036106 | 0.689938   |
| 2^16 | 8.785998 | 0.622057   | 9.658737 | 0.657884   | 9.727406 | 0.691300   |
| 2^17 | 9.413363 | 0.627365   | 10.31735 | 0.658613   | 10.41712 | 0.689713   |
| 2^18 | 10.04325 | 0.629887   | 10.97757 | 0.660217   | 11.10571 | 0.688590   |
| 2^19 | 10.67798 | 0.634735   | 11.63786 | 0.660288   | 11.79362 | 0.687912   |
| 2^20 | 11.31478 | 0.636795   | 12.30029 | 0.662431   | 12.4794  | 0.685775   |
| 2^21 | 11.95511 | 0.640335   | 12.96262 | 0.662330   | 13.16599 | 0.686594   |
| 2^22 | 12.59794 | 0.642821   | 13.62641 | 0.663791   | 13.85091 | 0.684920   |
| 2^23 | 13.24310 | 0.645163   | 14.29092 | 0.664509   | 14.53591 | 0.685001   |

Table 2 continued.

|  | 4 | difference | 10 | difference | 13 | difference |
|---|---|---|---|---|---|---|
| 2^0 |  |  |  |  |  |  |
| 2^1 |  |  |  |  |  |  |
| 2^2 |  |  |  |  |  |  |
| 2^3 |  |  |  |  |  |  |
| 2^4 | 0 |  |  |  |  |  |
| 2^5 | 0.693147 |  |  |  |  |  |
| 2^6 | 1.945910 | 1.252763 |  |  |  |  |
| 2^7 | 2.639057 | 0.693147 |  |  |  |  |
| 2^8 | 3.526361 | 0.887303 |  |  |  |  |
| 2^9 | 4.262680 | 0.736319 |  |  |  |  |
| 2^10 | 5.023881 | 0.761201 | 0 |  |  |  |
| 2^11 | 5.783825 | 0.759945 | 0.693147 |  |  |  |
| 2^12 | 6.505784 | 0.721959 | 1.945910 | 1.252763 |  |  |
| 2^13 | 7.247793 | 0.742009 | 2.708050 | 0.762140 | 0 |  |
| 2^14 | 7.960673 | 0.712880 | 3.610918 | 0.902868 | 0.693147 |  |
| 2^15 | 8.687948 | 0.727276 | 4.430817 | 0.819899 | 1.945910 | 1.252763 |
| 2^16 | 9.404179 | 0.716231 | 5.231109 | 0.800292 | 2.708050 | 0.762140 |
| 2^17 | 10.11787 | 0.713694 | 6.040255 | 0.809146 | 3.610918 | 0.902868 |
| 2^18 | 10.82862 | 0.710746 | 6.810142 | 0.769888 | 4.430817 | 0.819899 |
| 2^19 | 11.53721 | 0.708589 | 7.590347 | 0.780204 | 5.231109 | 0.800292 |
| 2^20 | 12.24503 | 0.707821 | 8.348775 | 0.758428 | 6.042633 | 0.811524 |
| 2^21 | 12.94923 | 0.704197 | 9.108086 | 0.759312 | 6.817831 | 0.775198 |
| 2^22 | 13.6536 | 0.704375 | 9.864227 | 0.756141 | 7.599902 | 0.782071 |
| 2^23 | 14.3556 | 0.702004 | 10.61325 | 0.749019 | 8.360305 | 0.760403 |

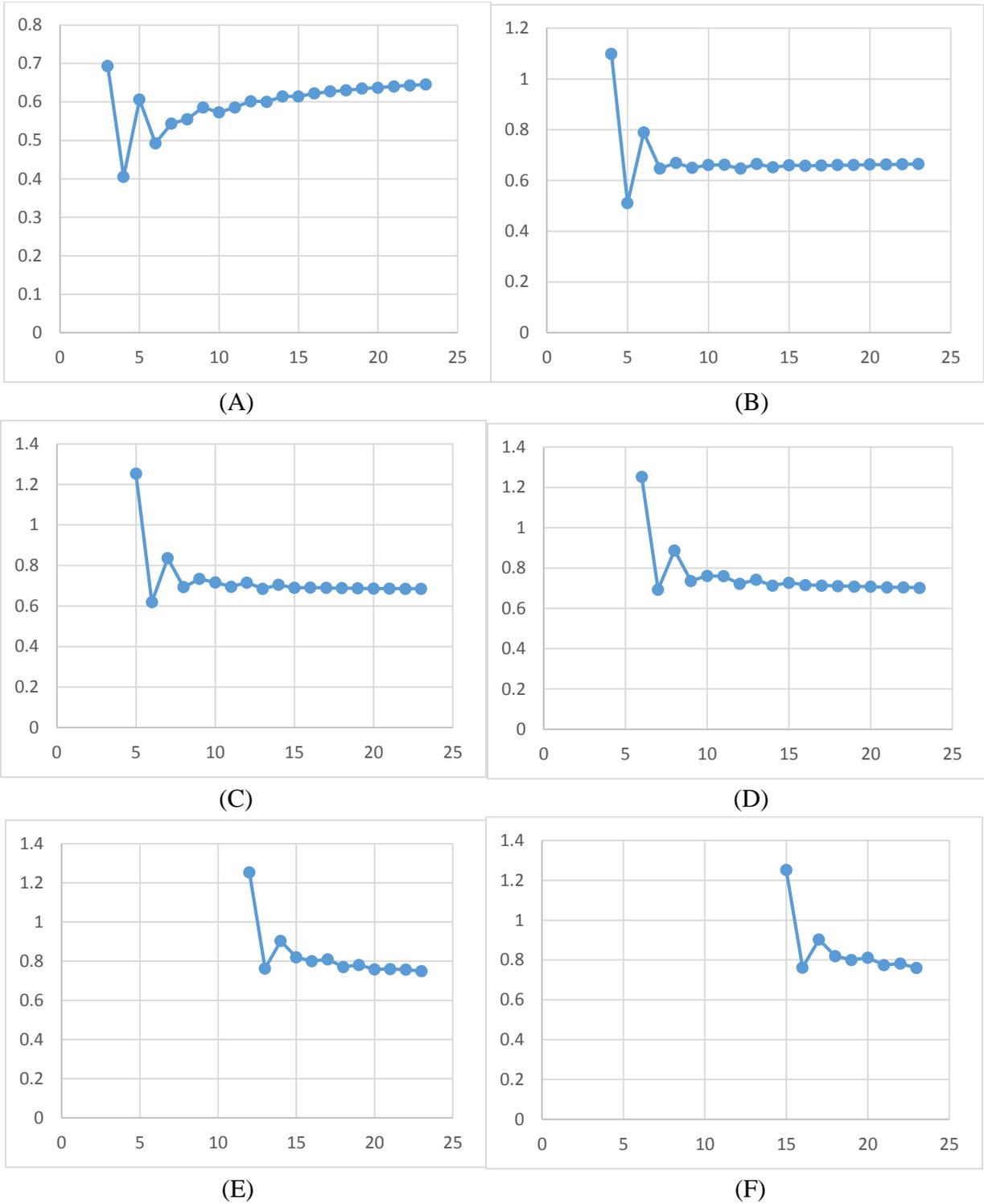

Fig 7. Difference of neighbor Y values of prime (A), 2D (B), 3D (C), 4D (D), 10D (E), and 13D column (F) in Table 2.

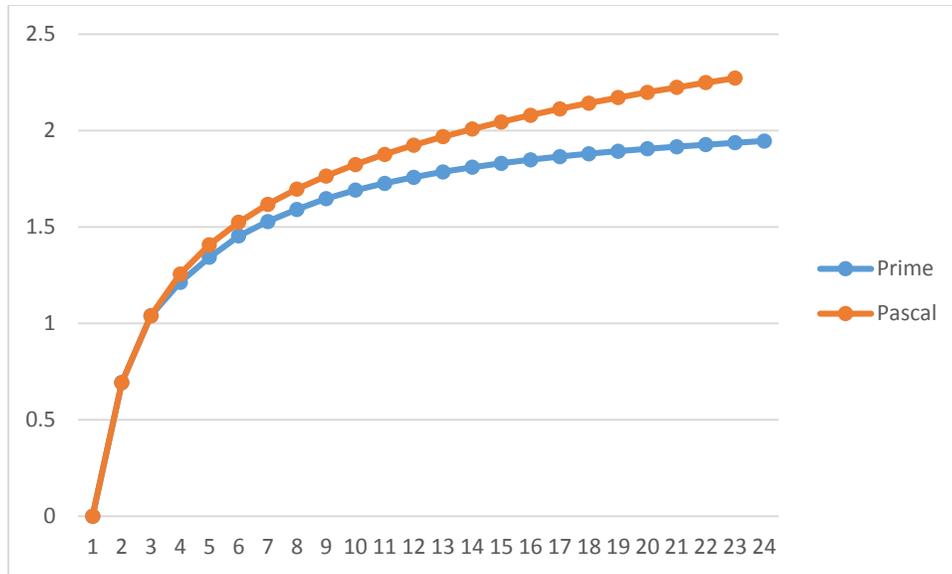

Fig 8. The entropy of the first 24 rows of dimensions distribution triangle and the first 23 rows of Pascal's triangle

Fig 8 shows that the entropy of dimension distributions is smaller than the distribution of row in Pascal's triangle and increases slower in each $2^n$ ($n \geq 4$) space. Comparing to Pascal's triangle, the dimension distributions are more concentrated and less random.

Question 4: Is the entropy increase of dimensions distribution triangle converging?

## 7. Summary

A triangle is constructed by counting all the natural numbers via their number of prime factors. The number of prime factors are defined as the 'dimension' of natural numbers. The natural number range [1, $n$] is defined to be '$n$ space'. Then the distributions of the dimensions in each $2^n$ ($n \geq 1$, $n \in \mathbb{N}$) space assemble the dimensions distribution triangle.

The triangle shows several interesting phenomena: 1) For each row, the last several terms are fixed; 2) the logarithm of the distributions looks quite disciplinary; 3) For $n \leq 23$, the maximum column shifts from the prime column to the 2D column and then the 3D column; 4) the entropy of distributions is smaller and increase slower than the distributions in Pascal's triangle. The entropy increase looks like converging.

Some questions raised from the superficial analysis of above phenomena, the answers of these questions may relate to the solution of prime number theorem. Due to the very limited mathematical knowledge of author, further explorations, e.g. connect the dimension distributions to the offset logarithmic integral function Li(x)[2], of this triangle have not been done. It will be great if someone could figure out something more of the triangle.

**Reference**

[1] Pascal's triangle. *In Wikipedia*. Retrieved Jan 1, 2014, from
   http://en.wikipedia.org/wiki/Pascal%27s_triangle
[2] Prime number theorem. *In Wikipedia*. Retrieved Jan 1, 2014, from
   http://en.wikipedia.org/wiki/Prime_number_theorem
[3] Entropy (statistical thermodynamics). *In Wikipedia*. Retrieved Jan 1, 2014, from
   http://en.wikipedia.org/wiki/Entropy_(statistical_thermodynamics)